\documentclass[a4paper,12pt,english]{article}
\usepackage{babel}
\usepackage{amsthm,amsmath}
\usepackage{amssymb}
\usepackage[all]{xypic}
\usepackage{enumerate}
\usepackage{a4wide}
\usepackage{hyperref}
\usepackage{latexsym}
\usepackage{enumerate}

\newtheorem{theorem}{Theorem}[section]
\newtheorem{proposition}[theorem]{Proposition}
\newtheorem{corollary}[theorem]{Corollary}
\newtheorem{lemma}[theorem]{Lemma}
\newtheorem*{theorem*}{Theorem}
\newtheorem*{proposition*}{Proposition}
\newtheorem*{corollary*}{Corollary}
\newtheorem*{lemma*}{Lemma}
\theoremstyle{definition}

\newtheorem{remark}[theorem]{Remark}
\newtheorem*{remark*}{Remark}

\newtheorem*{definition*}{Definition}

\newcommand{\cat}[1]{\mathcal{#1}}
\newcommand{\coring}[1]{\mathfrak{#1}}
\newcommand{\tensor}[1]{\otimes_{#1}}
\newcommand{\rcomod}[1]{\mathsf{Comod}_{#1}}
\newcommand{\rmod}[1]{\mathsf{Mod}_{#1}}

\newcommand{\rmodu}[1]{\overline{\mathsf{Mod}}_{#1}}

\newcommand{\cotensor}[1]{\square_{#1}}

\renewcommand{\hom}[3]{\mathrm{Hom}_{#1}(#2,#3)}

\newcommand{\rend}[2]{\mathrm{End}({#2}_{#1})}
\newcommand{\lend}[2]{\mathrm{End}({}_{#1}#2)}

\newcommand{\dostensor}[3]{#1 \tensor{#2} #3}
\newcommand{\trestensor}[5]{#1 \tensor{#2} #3 \tensor{#4} #5}
\newcommand{\fourtensor}[7]{#1 \tensor{#2} #3 \tensor{#4} #5 \tensor{#6} #7}
\newcommand{\abrir}[1]{e_{#1}\tensor{A}e^*_{#1}}
\newcommand{\coalg}[2]{{#1}_{#2}}

\begin{document}
\title{Comonads and Galois corings\footnote{Supported by the research project ``Algebraic Methods in Non Commutative Geometry'',
with financial support of the grant MTM2004-01406 from the DGICYT and FEDER}}
\author{J. G\'omez-Torrecillas \\
\normalsize Departamento de \'{A}lgebra \\ \normalsize Universidad
de Granada\\ \normalsize E18071 Granada, Spain \\
\normalsize e-mail: \textsf{gomezj@ugr.es} }

\date{}
\maketitle

\section*{Introduction}

The notion of a coring was introduced by M. E. Sweedler in
\cite{Sweedler:1975} with the objective of formulating and proving a
predual to the Jacobson-Bourbaki theorem for extensions of division
rings. A fundamental argument in \cite{Sweedler:1975} is the
following: given division rings $E \subseteq A$, each coideal $J$ of
the $A$--coring $A \tensor{E} A$ gives rise to a factor coring
$\coring{C} = A\tensor{E} A / J$. If $g \in \coring{C}$ denotes the
group-like element $1 \tensor{E} 1 + J$, then $D = \{ a \in A : ag =
ga \}$ is an intermediate division ring $E \subseteq D \subseteq A$.
Moreover, we have the canonical homomorphism of $A$--corings $\zeta
: A \tensor{D} A \rightarrow \coring{C}$ which sends $1 \tensor{D}
1$ onto $g$. It follows from \cite[2.2 Fundamental
Lemma]{Sweedler:1975} that $\zeta$ is an isomorphism of
$A$--corings. This fact is basic to establish the bijective
correspondence between coideals of $A \tensor{E} A$ and intermediate
extensions $E \subseteq D \subseteq A$ \cite[Fundamental
Theorem]{Sweedler:1975}. Sweedler's fundamental Lemma just referred
can be replaced by the fact that the coring $A \tensor{D} A$ turns
out to be simple cosemisimple \cite[Theorem
4.4]{ElKaoutit/Gomez/Lobillo:2004}, \cite[Theorem 3.2, Theorem
4.3]{ElKaoutit/Gomez:2003}, \cite[28.21]{Brzezinski/Wisbauer:2003}
or, alternatively, that $A$ is, as a right $\coring{C}$--comodule, a
simple generator of the category of all right
$\coring{C}$--comodules. Thus, we see that what is behind \cite[2.1
Fundamental Theorem]{Sweedler:1975} can be expressed in categorical
terms. In fact, this idea has been recently exploited to state a
generalization of Sweedler's theory for simple artinian rings
\cite{Cuadra/Gomez:2005arXiv}. In this work we show that the idea of
obtaining an isomorphism of corings from categorical properties can
be ultimately formulated in terms of comonads.

Each group-like element $g$ of a coring $\coring{C}$ over a ring
with unit $A$ gives \cite{Brzezinski:2002} a canonical homomorphism
of $A$--corings $\mathbf{can} : A \tensor{B} A \rightarrow
\coring{C}$, which sends $1 \tensor{B} 1$ onto $g$, where $B$ is the
subring of $A$ consisting of all $g$--coinvariant elements, and
$A\tensor{B} A$ is the Sweedler canonical coring
\cite{Sweedler:1975}. The group-like element $g$ also provides a
pair of adjoint functors between the category $\rmod{B}$ of right
$B$--modules and the category $\rcomod{\coring{C}}$ of right
$\coring{C}$--comodules \cite{Brzezinski:2002}. The left adjoint is
defined as a tensor product functor $- \tensor{B} A : \rmod{B}
\rightarrow \rcomod{\coring{C}}$, using the right
$\coring{C}$--comodule structure defined on $A$ by $g$. The right
adjoint may be interpreted as the functor $\hom{\coring{C}}{A}{-} :
\rcomod{\coring{C}} \rightarrow \rmod{B}$. T. Brzezi\'nski proves
\cite[Theorem 5.6]{Brzezinski:2002} that, with $\coring{C}$ flat as
a left $A$--module, this adjunction is an equivalence of categories
if, and only if, $\coring{C}$ is Galois (i.e., $\mathbf{can}$ is an
isomorphism) and $A$ is faithfully flat as a left $B$--module. The
canonical map $\mathbf{can}$ can be interpreted as a comonad in the
following way: The $A$--coring $\coring{C}$ gives rise to a comonad
on $\rmod{A}$ built over the functor $- \tensor{A} \coring{C}$
\cite[18.28]{Brzezinski/Wisbauer:2003}. On the other hand, the
adjunction associated to the ring extension $B \subseteq A$
determines \cite[Section 3.1]{Barr/Wells:2002} another comonad on
$\rmod{A}$. In this way, the canonical map $\mathbf{can}$ leads to
an homomorphism of comonads $- \tensor{A} \mathbf{can} : -
\tensor{A} A \tensor{B} A \rightarrow - \tensor{A} \coring{C}$. If
$(\coring{C},g)$ is a Galois coring, then these comonads are
isomorphic and the functor $- \tensor{B} A : \rmod{B} \rightarrow
\rcomod{\coring{C}}$ is, up to natural isomorphisms, the
Eilenberg-Moore comparison functor \cite[Section
3.2]{Barr/Wells:2002}. Thus, one of the implications of
\cite[Theorem 5.6]{Brzezinski:2002} may be obtained as a consequence
of Beck's Theorem \cite[Section 3.3]{Barr/Wells:2002}. This seems
not to be the case of the reciprocal implication. In fact, to deduce
that the canonical map $\mathbf{can}$ is an isomorphism from the
fact that $- \tensor{B} A : \rmod{B} \rightarrow
\rcomod{\coring{C}}$ is an equivalence one needs an independent
argument, which rests on precise relationship between the canonical
map and the counit of the adjunction given by the functors  $ -
\tensor{B} A$  and $\hom{\coring{C}}{A}{-}$. More precisely, to
deduce that $\coring{C}$ is Galois, it is enough to assume that the
counit is an isomorphism \cite[18.26]{Brzezinski/Wisbauer:2003},
\cite[Lemma 3.1]{ElKaoutit/Gomez:2003}, arising the condition of
being Galois as part of a characterization of the full and faithful
character of the functor $\hom{\coring{C}}{A}{-}$
\cite[18.27]{Brzezinski/Wisbauer:2003}, \cite[Theorem
3.8]{Caenepeel/DeGroot/Vercruysse:unp2005}, \cite[Remark
3.7]{ElKaoutit/Gomez:2003}. In fact, the referred results are proved
in the more general framework of the comatrix corings, introduced in
\cite{ElKaoutit/Gomez:2003}, where the role of the group-like
element $g$ is played by a right $\coring{C}$--comodule $\Sigma$
which is finitely generated and projective as a right module over
$A$, and $B = \rend{\coring{C}}{\Sigma}$. The generalization of
\cite[Theorem 5.6]{Brzezinski:2002} in this framework was proved in
\cite[Theorem 3.2]{ElKaoutit/Gomez:2003}. Our aim in this paper is
to investigate which aspects of these results admit a formulation in
terms of comonads. This approach hopefully will permit of focusing
in what is specific in each particular future situation, having some
relevant general results for granted.

Starting from a comonad $G$ on a category $\cat{A}$, and a functor
$L : \cat{B} \rightarrow \cat{A}$ with a right adjoint $R : \cat{A}
\rightarrow \cat{B}$, there is a bijective correspondence between
functors $K : \cat{B} \rightarrow \coalg{\cat{A}}{G}$ that factorize
throughout $L$  and homomorphisms of comonads $\varphi : LR
\rightarrow G$ \cite[Theorem II.1.1]{Dubuc:1970}. The notation
$K_{\varphi}$ will be used to refer to this dependence. Moreover, it
can be deduced from \cite[Theorem A.1]{Dubuc:1970} that $K_{\varphi}
: \cat{B} \rightarrow \coalg{\cat{A}}{G}$ admits a right adjoint
$D_{\varphi} : \coalg{\cat{A}}{G} \rightarrow \cat{B}$ under mild
conditions. We indicate elementary proofs of these facts, with the
aim of making them more accessible to non specialists in Category
Theory (see Proposition \ref{correspondencia}, Theorem \ref{Kphi}
and Proposition \ref{adjuncion}). We will prove that $D_{\varphi}$
is full and faithful if and only if $\varphi$ is an isomorphism and
$L$ preserves some equalizers (Theorem \ref{Descent}), and we will
conclude our general results by characterizing when $K_{\varphi}$
establishes an equivalence between the categories $\cat{B}$ and
$\coalg{\cat{A}}{G}$ (Theorem \ref{comonadic}). Obviously, the
functors characterized in this way are, a fortiori, comonadic but,
in contrast with the approach of Beck's Theorem, the comonad $G$ is
here given beforehand, and each functor $K_{\varphi}$ corresponds to
a ``representation'' of $G$. Beck's theorem deals with the situation
where $\varphi$ is the identity, that is, $G = LR$ is the comonad
associated to the adjunction.

Our point of view is motivated in part by the situation where an
entwining structure between an algebra and a coalgebra is given,
together with an entwined module \cite{Brzezinski/Majid:1998}. The
comonad $G$ is then given by the coring associated to the entwining
structure \cite[Section 2]{Brzezinski:2002}, and the functor
$K_{\varphi}$ is defined by an entwined module, which is noting but
a comodule over the aforementioned coring. It seems natural to study
the relationship between the category of entwined modules and the
category of modules over the subring of coinvariants defined by the
entwined module and, in particular, what the structure of the
comonad $G$ is when these categories are equivalent.

We will apply our general theory to the study of comodules over
corings over firm rings. This illustrates how the results on
comonads given in the Section \ref{candefinido} significantly
simplify the treatment of some relevant aspects of the comatrix
corings and the Galois comodules investigated in
\cite{ElKaoutit/Gomez:2003}, \cite{Gomez/Vercruysse:2005arXiv} or
\cite{Wisbauer:2004arXiv}. We show how the notion of Galois comodule
without finiteness conditions as introduced in
\cite{Wisbauer:2004arXiv} fits perfectly in the general categorical
setting. We derive, in particular, some new results on full and
faithful functors and equivalences between categories of modules and
comodules given by Galois comodules without assuming a priori
finiteness conditions (Theorem \ref{fielmenteplano}, Theorem
\ref{GE}).

\section{Functors with values in coalgebras and homomorphisms of comonads}\label{candefinido}

Let $(G,\Delta, \varepsilon)$ be a comonad (or cotriple) on a
category $\cat{A}$, that is, a functor $G : \cat{A} \rightarrow
\cat{A}$ together with two natural transformations $\Delta : G
\rightarrow G^2$ and $\varepsilon : G \rightarrow id_{\cat{A}}$
such that the diagrams
\begin{equation*}
\xymatrix{G \ar^{\Delta}[r] \ar_{\Delta}[d] & G^2 \ar^{G\Delta}[d]
\\ G^2 \ar_{\Delta G}[r] & G^3 } \qquad \xymatrix{G & G^2 \ar_{\varepsilon G}[l] \ar^{G \varepsilon}[r] & G \\
& \ar@{=}[ul] G \ar|{\Delta}[u] \ar@{=}[ur] & }
\end{equation*}
are commutative \cite{Eilenberg/Moore:1965}, \cite[Chapter
3]{Barr/Wells:2002}. We will follow as much as possible
\cite{Barr/Wells:2002}, understanding each statement on monads
(triples) automatically as its version for comonads. Consider a
functor $L : \cat{B} \rightarrow \cat{A}$ with a right adjoint $R
: \cat{A} \rightarrow \cat{B}$. If  $\eta : id_{\cat{B}}
\rightarrow RL$ is the unit of the adjunction, and $\epsilon : LR
\rightarrow id_{\cat{A}}$ is its counit, then we have
\cite[Proposition 3.1.2]{Barr/Wells:2002} the comonad $(LR,
\delta, \epsilon)$ on $\cat{A}$, where $\delta = L\eta R$. Let us
recall the category $\coalg{\cat{A}}{G}$ of $G$--coalgebras
\cite[Section 3.1]{Barr/Wells:2002}, whose objects are pairs
$(X,x)$ consisting of an object $X$ of $\cat{A}$ and a morphism $x
: X \rightarrow GX$ such that
\begin{equation*}
Gx \circ x = \Delta_X \circ x, \qquad \varepsilon_X \circ x = id_X
\end{equation*}
Given $G$--coalgebras $(X,x), (X',x')$, the morphisms $f : X
\rightarrow X'$ in $\cat{A}$ such that $Gf \circ x = x' \circ f$
form the set of homomorphisms $\hom{\coalg{\cat{A}}{G}}{X}{X'}$ in
$\coalg{\cat{A}}{G}$ from $(X,x)$ to $(X',x')$.

Let $U : \coalg{\cat{A}}{G} \rightarrow \cat{A}$ denote the
forgetful functor. We will consider those functors $K : \cat{B}
\rightarrow \coalg{\cat{A}}{G}$ such that the diagram
\begin{equation}\label{UKL}
\xymatrix{\cat{B} \ar^{K}[r] \ar_{L}[dr] & \coalg{\cat{A}}{G}
\ar^{U}[d] \\
& \cat{A}}
\end{equation}
is commutative. We are specially interested in the case where $K$
provides an equivalence of categories between $\cat{B}$ and
$\coalg{\cat{A}}{G}$. We start by giving an elementary proof of the
bijective correspondence, formulated in \cite[Theorem
II.1.1]{Dubuc:1970}, between the functors $K$ that make commute the
diagram \eqref{UKL} and the homomorphisms of comonads $\varphi : LR
\rightarrow G$. This correspondence will be a consequence of the
following Proposition \ref{correspondencia} which, in addition,
establish some technical facts that will be needed later. Recall
\cite[Section 3.6]{Barr/Wells:2002} that a homomorphism of comonads
from $LR$ to $G$ is a natural transformation $\varphi : LR
\rightarrow G$ such that $\Delta \varphi = \varphi^2 \delta$ and
$\varepsilon \varphi = \epsilon$. The one-to-one correspondence
between the mathematical objects described in the statements (A) and
(C) of Proposition \ref{correspondencia} can be deduced from
\cite[Proposition II.1.4]{Dubuc:1970}.

\begin{proposition}\label{correspondencia}
There exist bijective correspondences between: \\

\noindent (A) Homomorphisms of comonads from $(LR, L\eta R,
\epsilon)$ to
$(G,\Delta,\varepsilon)$. \\

\medskip

 \noindent (B) Natural transformations $\xymatrix{R
\ar^{\alpha}[r] & RG}$ such that the following diagrams commute:
\begin{equation}\label{comod}
\xymatrix{R \ar^{\alpha}[r] \ar^{\alpha}[d] & RG \ar^{R\Delta}[d] \\
RG \ar^{\alpha G}[r] & RG^2 } \qquad \xymatrix{ R \ar^{\alpha}[r]
\ar@{=}[dr] & RG \ar^{R \varepsilon}[d] \\
& R }
\end{equation}
and

\noindent (C)  Natural transformations $\xymatrix{L \ar^{\beta}[r]
& GL}$ such that the following diagrams commute:
\begin{equation}\label{comodb}
\xymatrix{L \ar^{\beta}[r] \ar^{\beta}[d] & GL \ar^{\Delta L}[d] \\
GL \ar^{ G \beta }[r] & G^2L } \qquad \xymatrix{ L \ar^{\beta}[r]
\ar@{=}[dr] & GL \ar^{\varepsilon L}[d] \\
& L }
\end{equation}
\end{proposition}
\begin{proof}
We first prove the one-to-one correspondence between the natural
transformations described in $(A)$ and $(B)$. Given a homomorphism
of comonads $\varphi : LR \rightarrow G$, define the natural
transformation
\begin{equation}\label{rhophi}
\xymatrix{ R \ar@/_1pc/_{\alpha}[rr] \ar^{\eta R}[r] & RLR \ar^{R
\varphi}[r] & RG}
\end{equation}
To check that the first diagram in \eqref{comod} commutes,
consider an object $X$ of $\cat{A}$ and compute:
\begin{equation*}
\begin{array}{lclr}
R\Delta_X \circ \alpha_X & = & R\Delta_X \circ R\varphi_X \circ
\eta_{RX} & \\
& = & R\varphi^2_X \circ RL\eta_{RX} \circ \eta_{RX} & (\varphi
\hbox{ is of comonads}) \\
 & = & R\varphi_{GX} \circ RLR\varphi_X \circ RL\eta_{RX} \circ
 \eta_{RX} & (\varphi^2_X = \varphi_{GX} \circ LR\varphi_X) \\
 & = & R\varphi_{GX} \circ RLR\varphi_X \circ \eta_{RLRX} \circ
 \eta_{RX} & (\eta \hbox{ is natural}) \\
 & = & R\varphi_{GX} \circ \eta_{RGX} \circ R\varphi_X \circ \eta_{RX} &
 (\eta \hbox{ is natural}) \\
 & = & \alpha_{GX} \circ \alpha_X
\end{array}
\end{equation*}
For the second diagram in \eqref{comod}, we have:
\begin{equation*}
\begin{array}{lclr}
R\varepsilon_X \circ \alpha_X & = & R\varepsilon_X \circ
R\varphi_X
\circ \eta_{RX} & \\
& = & R\epsilon_X \circ \eta_{RX} & (\varphi \hbox{ is of
comonads}) \\
& = & id_{RX} & (\hbox{by the adjunction})
\end{array}
\end{equation*}
Conversely, given a natural transformation $\alpha : R \rightarrow
RG$ satisfying \eqref{comod}, we define the natural transformation
\begin{equation}\label{can}
\xymatrix{LR \ar^{L\alpha}[r] \ar@/_1pc/_{\varphi}[rr] & LRG
\ar^{\epsilon G}[r] & G}
\end{equation}
To show that $\varphi$, defined in \eqref{can}, is  a homomorphism
of comonads, we need to use that, for each object $X$ of
$\cat{A}$, one has the identities
\begin{equation}\label{can2}
\varphi^2_X = \epsilon_{G^2X} \circ L\alpha_{GX} \circ
LR\epsilon_{GX} \circ LRL\alpha_X
\end{equation}
and
\begin{equation*}
\varphi^2_X = G\epsilon_{GX} \circ GL\alpha_X \circ
\epsilon_{GLRX} \circ L\alpha_{RLX},
\end{equation*}
by definition of $\varphi^2$. Do the following computation:
\begin{equation*}
\begin{array}{lclr}
\varphi^2_X \circ L\eta_{RX} & = & \epsilon_{G^2X} \circ
L\alpha_{GX} \circ LR\epsilon_{GX} \circ LRL\alpha_X \circ
L\eta_{RX}
& (\hbox{by \eqref{can2}}) \\
& = &  \epsilon_{G^2X} \circ L\alpha_{GX} \circ LR\epsilon_{GX}
\circ L\eta_{RGX} \circ L\alpha_X & (\eta \hbox{ is natural}) \\
& = & \epsilon_{G^2X} \circ L\alpha_{GX}  \circ L\alpha_X &
(R\epsilon_{GX} \circ \eta_{RGX} = id_{RGX}) \\
& = & \epsilon_{G^2X} \circ LR\Delta_X \circ L\alpha_X &
(\alpha_{GX} \circ \alpha_X = R\Delta_X \circ \alpha_X) \\
& = & \Delta_X \circ \epsilon_{GX} \circ L\alpha_X & (\epsilon
\hbox{ is natural}) \\
& = & \Delta_X \circ \varphi_X
\end{array}
\end{equation*}
This proves that the first condition for $\varphi$ to be a
homomorphism of comonads holds. For the second condition, we compute
as follows:
\begin{equation*}
\begin{array}{lclr}
\varepsilon_X \circ \varphi_X & = & \varepsilon_X \circ
\epsilon_{GX} \circ L\alpha_X & \\
& = & \epsilon_X \circ LR\varepsilon_X \circ L\alpha_X & (\epsilon
\hbox{ is natural}) \\
& = & \epsilon_X & (R\varepsilon_X \circ \alpha_X = id_{RX})
\end{array}
\end{equation*}
Now, let  $\varphi : LR \rightarrow G$ be a homomorphism of
comonads and consider the natural transformation $\alpha : R
\rightarrow RG$ defined in \eqref{rhophi}. If we consider the
homomorphism of comonads, say $\varphi'$, defined, from this
$\alpha$, in \eqref{can}, then we will see that it coincides with
the original $\varphi$. For this, we compute, for an object $X$ of
$\cat{A}$:
\begin{equation*}
\begin{array}{lclr}
\varphi'_X & = & \epsilon_{GX} \circ L\alpha_X & \\
& = & \epsilon_{GX} \circ LR\varphi_X \circ L\eta_{RX} & \\
& = & \varphi_X \circ \epsilon_{LRX} \circ L\eta_{RX} & (\epsilon
\hbox{ is natural}) \\
& = & \varphi_X &
\end{array}
\end{equation*}
Finally, we have to see that, given a natural transformation $\alpha
: R \rightarrow RG$ subject to \eqref{comod}, and defining first
$\varphi$ according to \eqref{can}, and the new natural
transformation, say $\alpha'$, by \eqref{rhophi}, we reobtain the
original $\alpha$. This follows from the following computation, for
$X$ an object of $\cat{A}$:
\begin{equation*}
\begin{array}{lclr}
\alpha'_X & = & R\varphi_X \circ \eta_{RX} & \\
& = & R\epsilon_{GX} \circ RL\alpha_X \circ \eta_{RX} & (\eta
\hbox{
is natural}) \\
& = & R\epsilon_{GX} \circ \eta_{RGX} \circ \alpha_X & \\
& = & \alpha_X &
\end{array}
\end{equation*}
The correspondence between the natural transformations described
in $(A)$ and $(C)$ can be proved in a similar way (details can be
found in \cite[Proposici\'{o}n 1.1]{Gomez:2006arXiv}). The assignment
goes as follows: starting from a homomorphism of comonads
$\varphi: LR \rightarrow G$, we define the natural transformation
\begin{equation}\label{lambdaphi}
\xymatrix{ L \ar@/_1pc/_{\beta}[rr] \ar^{L \eta}[r] & LRL \ar^{
\varphi L}[r] & GL}
\end{equation}
which turns out to make commute the diagrams \eqref{comodb}.
Conversely, to each natural transformation $\beta : L \rightarrow
GL$ making commute the diagrams in \eqref{comodb}, we assign the
homomorphism of comonads $\varphi$ given by
\begin{equation}\label{canr}
\xymatrix{LR \ar^{\beta R}[r] \ar@/_1pc/_{\varphi}[rr] & GLR \ar^{
G \epsilon}[r] & G}
\end{equation}
\end{proof}

Given a homomorphism of comonads $\varphi : LR \rightarrow G$, and
the corresponding natural transformation $\beta$ given by
\eqref{lambdaphi}, we can define the functor $K_{\varphi} :
\cat{B} \rightarrow \coalg{\cat{A}}{G}$ defined as
\begin{equation*}
K_{\varphi}Y = (LY,\beta_Y) = (LY,\varphi_{LY} \circ L\eta_Y)
\end{equation*}
on objects $Y$ of $\cat{B}$, and as $K_{\varphi}f = Lf$ on
morphisms $f$ of $\cat{B}$. The diagrams \eqref{comodb} show that
$(LY,\beta_Y)$ is a $G$--coalgebra and the naturalness of $\beta$
implies that $Lf$ is a homomorphism of $G$--coalgebras. This
functor $K_{\varphi}$ satisfies that $U K_{\varphi} = L$.
Conversely, given a functor $K : \cat{B} \rightarrow
\coalg{\cat{A}}{G}$ making commute the diagram \eqref{UKL},
define, for each object $Y$ of $\cat{B}$, $\beta_Y : LY
\rightarrow GLY$ as the structure morphism of the $G$--coalgebra
$KY$. It is easy to check that this assignment defines a natural
transformation $\beta : L \rightarrow GL$ which makes commute
\eqref{comodb}. We thus deduce from Proposition
\ref{correspondencia}:

\begin{theorem}\label{Kphi}\cite[Theorem II.1.1]{Dubuc:1970}.
Given a comonad $G$ on a category $\cat{A}$, and a functor $L :
\cat{B} \rightarrow \cat{A}$, if $L$ has a right adjoint $R :
\cat{A} \rightarrow \cat{B}$, then there exists a bijective
correspondence between functors $K : \cat{B} \rightarrow
\coalg{\cat{A}}{G}$ making commute \eqref{UKL} and homomorphisms of
comonads $\varphi : LR \rightarrow G$.
\end{theorem}

The situation described in Proposition \ref{correspondencia}
requires a suitable terminology. Thus, given a homomorphism of
comonads $\varphi : LR \rightarrow G$, we will refer to the natural
transformations $\alpha : R \rightarrow RG$ and $\beta : L
\rightarrow GL$ as the \emph{co-induced representation} (resp.
\emph{induced representation}) of $\varphi$. When the starting datum
is either a natural transformation $\alpha$ subject to the
conditions \eqref{comod}, or a natural transformation $\beta$ making
commute the diagrams in \eqref{comodb}, then the corresponding
homomorphism of comonads $\varphi : LR \rightarrow G$ will be called
\emph{canonical transformation} associated to $\alpha$ (resp.
$\beta$).

We next study when $K_{\varphi}$ has a right adjoint. The
following Proposition \ref{adjuncion} can be deduced from
\cite[Theorem A.1]{Dubuc:1970}. We give an elementary proof in our
case.

\begin{proposition}\label{adjuncion}
Assume that for every $G$--coalgebra $(X,x)$ there exists in
$\cat{B}$ the equalizer of the pair of morphisms $\alpha_X, Rx :
RX \rightarrow RGX$. Then the functor $K_{\varphi} : \cat{B}
\rightarrow \coalg{\cat{A}}{G}$ has a right adjoint $D_{\varphi} :
\coalg{\cat{A}}{G} \rightarrow \cat{B}$, whose value at $(X,x)$ is
the equalizer
\begin{equation}\label{Dphi}
\xymatrix{D_{\varphi}X \ar^{eq_X}[rr] & & RX
\ar@<.5ex>^-{\alpha_X}[rr] \ar@<-.5ex>_-{Rx}[rr] & & RGX}
\end{equation}
\end{proposition}
\begin{proof}
Given objects $X$ of $\cat{A}$ and $Y$ of $\cat{B}$, denote by $\Phi
: \hom{\cat{A}}{LY}{X} \rightarrow \hom{\cat{B}}{Y}{RX}$ the
isomorphism of the adjunction, which, in terms of the counity
$\eta$, is defined by $\Phi (h) = Rh \circ \eta_Y$ for $h : LY
\rightarrow X$. Consider the following commutative diagram of maps
between sets:
\begin{equation}\label{adjdia}
\xymatrix{\hom{\coalg{\cat{A}}{G}}{K_{\varphi}Y}{X}
\ar@{-->}^{\Upsilon}[rr] \ar[d] & &
\hom{\cat{B}}{Y}{D_{\varphi}X} \ar[d] \\
\hom{\cat{A}}{LY}{X} \ar^{\Phi}[rr] \ar@<-1ex>_{x \circ - }[d]
\ar@<1ex>^{G(-) \circ \beta_Y}[d]& &
\hom{\cat{B}}{Y}{RX}  \ar@<-1ex>_{\hom{\cat{B}}{Y}{Rx}}[d] \ar@<1ex>^{\hom{\cat{B}}{Y}{\alpha_X}}[d]\\
\hom{\cat{A}}{LY}{GX} \ar^{\Phi}[rr]
 & &
\hom{\cat{B}}{Y}{RGX},}
\end{equation}
where $\Upsilon$ is given by the restriction of $\Phi$ to
$\hom{\coalg{\cat{A}}{G}}{K_{\varphi}Y}{X}$. In fact, the bottom
square commutes serially, in the sense of \cite[page
112]{Barr/Wells:2002} (a detailed verification can be found in
\cite[Proposici\'{o}n 2.1]{Gomez:2006arXiv}). Now, the vertical edges
are equalizers, the left one by definition of homomorphism of
$G$--coalgebras, and the right edge by the universal property of the
equalizer \eqref{Dphi}. Hence, the natural isomorphism $\Phi :
\hom{\cat{A}}{LY}{X} \rightarrow \hom{\cat{B}}{Y}{RX}$ induces, by
restriction, a natural isomorphism $\Upsilon :
\hom{\coalg{\cat{A}}{G}}{K_\varphi Y}{X} \rightarrow
\hom{\cat{B}}{Y}{D_{\varphi}X}$.
\end{proof}

\begin{remark}\label{unitcounit}
If we put in the diagram \eqref{adjdia} $X = K_{\varphi}Y$, for an
object $Y$ of $\cat{B}$, since $id_{K_{\varphi}Y}$ is a
homomorphism of $G$--coalgebras, we get that $\eta_Y =
\Phi(id_{K_{\varphi}Y})$ factorizes throughout
$D_{\varphi}K_{\varphi}Y$. In this way,  the unit of the
adjunction $K_{\varphi} \dashv D_{\varphi}$, denoted by
$\widehat{\eta}$, is uniquely determined at $Y$ by the universal
property of an equalizer, according to the following diagram:
\begin{equation}\label{unitphi} \xymatrix{D_{\varphi}K_{\varphi}Y \ar[r] &
RLY \ar@<.5ex>^-{\alpha_{LY}}[rr] \ar@<-.5ex>_-{R\beta_Y}[rr] & &
RGLY \\
& Y \ar@{-->}^{\widehat{\eta}_Y}[ul] \ar^{\eta_Y}[u] & & }
\end{equation}
To make explicit the counit $\widehat{\epsilon}$ at an object
$(X,x)$ of $\coalg{\cat{A}}{G}$, take $Y = D_{\varphi}X$ in
diagram \eqref{adjdia}. Then we obtain the morphism
\begin{equation*}
\widehat{\epsilon}_X = \Upsilon^{-1}(id_{D_{\varphi}X}) =
\Phi^{-1}(eq_X) = \epsilon_X \circ Leq_X
\end{equation*}
which is of $G$--coalgebras in view of \eqref{adjdia}. This
situation is resumed by the diagram in $\cat{A}$:
\begin{equation}\label{counidad1}
\xymatrix{K_{\varphi}D_{\varphi}X \ar^{Leq_X}[rr]
\ar@{-->}^{\widehat{\epsilon}_X}[drr] & & LRX
\ar@<.5ex>^-{L\alpha_X}[rr] \ar@<-.5ex>_-{LRx}[rr] \ar^{\epsilon_X}[d] & & LRGX \\
& & X, & & .}
\end{equation}
well understood that the underlying object in $\cat{A}$ of
$K_{\varphi}D_{\varphi}X$ is, by definition, $LD_{\varphi}X$.
\end{remark}

Our next aim is to prove that $D_{\varphi}$ is a full and faithful
functor if and only if $\varphi$ is an isomorphism of comonads and
the functor $L$ preserves the equalizers \eqref{Dphi}.

\begin{lemma}\label{clave3}
Assume that the equalizer \eqref{Dphi} exists for every
$G$--coalgebra $(X,x)$. Then the following diagram
\begin{equation}\label{clave3b}
\xymatrix{ LD_{\varphi}X \ar^{Leq_X}[rr]
\ar^{\widehat{\epsilon}_X}[d]& & LRX
\ar@<.5ex>^-{L\alpha_X}[rr] \ar@<-.5ex>_-{LRx}[rr] \ar^{\varphi_X}[d] & & LRGX \ar^{\varphi_{GX}}[d]\\
X \ar^{x}[rr] & & GX \ar@<.5ex>^-{\Delta_X}[rr]
\ar@<-.5ex>_-{Gx}[rr] & & G^2X  }
\end{equation}
is (serially) commutative.
\end{lemma}
\begin{proof}
That the right square commutes serially is an easy consequence of
the naturalness of $\varphi$ and of the fact that it is a
homomorphism of comonads. For the commutativity of the left
square, we do the following computation:
\begin{equation*}
\begin{array}{lclr}
x \circ \widehat{\epsilon}_X & = & x \circ \epsilon_X \circ
Leq_X & \\
& = & \epsilon_{GX} \circ LRx \circ Leq_X & (\epsilon \hbox{ is
natural}) \\
& = & \epsilon_{GX} \circ L\alpha_X  \circ Leq_X & (eq_X \hbox{ equalizes } (Rx,\alpha_X))\\
& = & \epsilon_{GX} \circ LR\varphi_X \circ L\eta_{RX} \circ
Leq_X & \\
& = & \varphi_X \circ \epsilon_{LRX} \circ L\eta_{RX} \circ
Leq_X & (\epsilon \hbox{ is natural}) \\
& = & \varphi_X \circ Leq_X &
\end{array}
\end{equation*}
\end{proof}

We are now ready to state the main result in this section. Recall
that a right adjoint functor is faithful and full if and only if
the counit is an isomorphism \cite[Proposition
3.4.1]{Borceux:1994}

\begin{theorem}\label{Descent}
Let $L : \cat{B} \rightarrow \cat{A}$ be a functor admitting a right
adjoint $R : \cat{A} \rightarrow \cat{B}$, and let $G : \cat{A}
\rightarrow \cat{A}$ be a comonad. Consider a functor $K_{\varphi} :
\cat{B} \rightarrow \coalg{\cat{A}}{G}$ that makes commute
\eqref{UKL} with corresponding homomorphism of comonads $\varphi :
LR \rightarrow G$, and let $\alpha : R \rightarrow RG$, $\beta : L
\rightarrow GL$ be its representations. Assume that for every
$G$--coalgebra $(X,x)$, there exists in $\cat{B}$ the equalizer of
$\alpha_X, Rx$. Then the right adjoint $D_{\varphi} :
\coalg{\cat{A}}{G} \rightarrow \cat{B}$ to the functor $K_{\varphi}$
is faithful and full if and only if $L$ preserves the equalizers of
the form \eqref{Dphi} and $\varphi$ is an isomorphism of comonads.
\end{theorem}
\begin{proof}
 Let $X$ be
any object of $\cat{A}$. An easy computation shows that
\begin{equation}\label{ecucontr}
\xymatrix{RX \ar^-{\alpha_X}[rr] & &
\ar@/^1pc/^-{R\varepsilon_X}[ll] RGX \ar@<.5ex>^-{\alpha_{GX}}[rr]
\ar@<-.5ex>_-{R\Delta_X}[rr] & & RG^2X
\ar@/^1.5pc/^-{RG\varepsilon_X}[ll]}
\end{equation}
is a contractible equalizer in the sense of \cite[Section
3.3]{Barr/Wells:2002}. The equalizer \eqref{ecucontr} shows that
$eq_{GX} = \alpha_X$ and $D_{\varphi}GX = RX$. From this, by
applying the functor $L$ to \eqref{ecucontr}, and taking
\eqref{counidad1} into account, we get that $\widehat{\epsilon}_{GX}
= \epsilon_{GX} \circ L\alpha_X = \varphi_X$. Thus, if $D_{\varphi}$
is full and faithful, then $\widehat{\epsilon}_{GX} = \varphi_X$ is
an isomorphism for every $G$--coalgebra $(X,x)$. Moreover, since the
bottom row in the commutative diagram \eqref{clave3b} is an
equalizer \cite[Proposition 3.3.4]{Barr/Wells:2002}, we get that its
top row is an equalizer as well. This means that $L$ preserves the
equalizers of the form \eqref{Dphi}. Conversely, if $\varphi$ is an
isomorphism and $L$ preserves all equalizers of the form
\eqref{Dphi}, then, by Lemma \ref{clave3}, $\widehat{\epsilon}_X$ is
an isomorphism for every $G$--coalgebra $(X,x)$, whence
$D_{\varphi}$ is full and faithful.
\end{proof}

Any isomorphism of comonads $\varphi : LR \rightarrow G$ induces
\cite[Theorem 3.3]{Barr/Wells:2002} an equivalence of categories
$\coalg{\cat{A}}{LR} \cong \coalg{\cat{A}}{G}$. By combining this
fact with Theorem \ref{Descent} and Beck's Theorem \cite[Theorem
3.10]{Barr/Wells:2002}, we may obtain Theorem \ref{comonadic}. We
prefer to include here an explicit proof.

\begin{theorem}\label{comonadic}
Let $L : \cat{B} \rightarrow \cat{A}$ be a functor with a right
adjoint $R : \cat{A} \rightarrow \cat{B}$, and $G : \cat{A}
\rightarrow \cat{A}$ any comonad. Consider a functor $K_{\varphi}
: \cat{B} \rightarrow \coalg{\cat{A}}{G}$ that makes commute
\eqref{UKL} with the corresponding homomorphism of comonads
$\varphi : LR \rightarrow G$, and let $\alpha : R \rightarrow RG$,
$\beta : L \rightarrow GL$ be its representations. Assume that for
every $G$--coalgebra $(X,x)$, there exists in $\cat{B}$ the
equalizer of $\alpha_X, Rx$. Then the functor $K_{\varphi}$ is an
equivalence of categories between $\cat{B}$ and
$\coalg{\cat{A}}{G}$ if and only if $L$ preserves the equalizers
of the form \eqref{Dphi}, reflects isomorphisms, and $\varphi$ is
an isomorphism of comonads.
\end{theorem}
\begin{proof}
We first observe that for every object $Y$ of $\cat{B}$, the unit
$\widehat{\eta}_Y$ of the adjunction $K_{\varphi} \dashv
D_{\varphi}$ stated in Proposition \ref{adjuncion} is given,
according to the Remark \ref{unitcounit}, by the equalizer in the
horizontal row in the commutative diagram
\begin{equation}\label{unitphi2}
\xymatrix{ & & & RLRLY \ar^-{R\varphi_{LY}}[dr] & \\
D_{\varphi}K_{\varphi}Y \ar[rr] & & RLY
\ar@<.5ex>^-{\alpha_{LY}}[rr] \ar@<-.5ex>_-{R\beta_Y}[rr]
\ar@<.5ex>^-{\eta_{RLY}}[ur] \ar@<-.5ex>_-{RL\eta_Y}[ur]& &
RGLY\\
& Y \ar^-{\widehat{\eta}_Y}[ul] \ar_{\eta_Y}[ur] & & &}
\end{equation}
If we apply the functor $L$ to the commutative diagram
\eqref{unitphi2} we obtain the diagram
\begin{equation}\label{Lunitphi2}
\xymatrix{ & & & LRLRLY \ar^-{LR\varphi_{LY}}[dr] \ar@/_2.5pc/_{\epsilon_{LRLY}}[ld]& \\
LD_{\varphi}K_{\varphi}Y \ar[rr] & & LRLY
\ar@/_1pc/_-{\epsilon_{LY}}[ld] \ar@<.5ex>^-{L\alpha_{LY}}[rr]
\ar@<-.5ex>_-{LR\beta_Y}[rr] \ar@<.5ex>^-{L\eta_{RLY}}[ur]
\ar@<-.5ex>_-{LRL\eta_Y}[ur] & &
LRGLY\\
& LY \ar^-{L\widehat{\eta}_Y}[ul] \ar_{L\eta_Y}[ur], & & &}
\end{equation}
commutative as well. Here, the morphisms $\epsilon_{LRLY}$ and
$\epsilon_{LY}$ make the diagonal row a contractible equalizer. If
$K_{\varphi}$ is an equivalence of categories, then its right
adjoint $D_{\varphi}$ is obviously faithful and full and, by
Theorem \ref{Descent}, $\varphi$ is a natural isomorphism and $L$
preserves the equalizers of the form \eqref{Dphi}. Since the
forgetful functor $U_G : \coalg{\cat{A}}{G} \rightarrow \cat{A}$
reflects isomorphisms \cite[Proposition 3.3.1]{Barr/Wells:2002},
we get from $L = U_G \circ K_{\varphi}$ that $L$ reflects
isomorphisms. Conversely, if $\varphi$ is a natural isomorphism
and $L$ preserves the equalizers \eqref{Dphi} and reflects
isomorphisms, then, by Theorem \ref{Descent}, the counit of the
adjunction $K_{\varphi} \dashv D_{\varphi}$ is an isomorphism.
From the diagram \eqref{Lunitphi2} we deduce that
$L\widehat{\eta}_Y$ is an isomorphism and, since $L$ reflects
isomorphisms, $\widehat{\eta}_Y$ must be an isomorphism. We have
thus proved that the unit of the adjunction $K_{\varphi} \dashv
D_{\varphi}$ is also a natural isomorphism. Therefore,
$K_{\varphi}$ is an equivalence of categories.
\end{proof}

\section{Corings over firm rings}

Let $A$ be a ring, which is not assumed to have a unit. By
$\rmodu{A}$ we denote the category of all right $A$--modules. We
may also consider left modules or bimodules over different rings.
The tensor product over $A$ will be denoted by $- \tensor{A} - $.
A right $A$--module $M$ is said to be \emph{firm}
\cite{Quillen:unp1997} if the ``multiplication'' map $\varpi_M^+ :
M \tensor{A} A \rightarrow M$ is bijective. Its inverse will be
denoted by $d_M^+ : M \rightarrow M \tensor{A} A$. Left firm
$A$--modules are defined analogously, with notations $\varpi_M^-$
and $d_M^-$ for the ``multiplication'' map and its inverse,
respectively. Obviously, $\varpi_A^+ = \varpi_A^-$, which implies,
in case of being bijective, that $d_A^+ = d_A^-$. We will say then
that $A$ is a \emph{firm ring}. If $A$ is firm, then the full
subcategory $\rmod{A}$ of $\rmodu{A}$ whose objects are all firm
modules is abelian \cite[(4.6)]{Quillen:unp1997}, \cite[Corollary
1.3]{Grandjean/Vitale:1998}, \cite[Proposition 2.7]{Marin:1998},
even thought, in general, that equalizers cannot be computed in
Abelian Groups, due essentially to the lack of exactness of the
functor $ - \tensor{A} A$.

Let $A$ be a firm ring. An \emph{$A$--coring} is a coalgebra in
the monoidal category of all firm $A$--bimodules (see, e.g.,
\cite[38.33]{Brzezinski/Wisbauer:2003}). Explicitly, an
$A$--coring is a firm $A$--bimodule $\coring{C}$ endowed with two
homomorphisms of $A$--bimodules $\Delta_{\coring{C}} : \coring{C}
\rightarrow \coring{C} \tensor{A} \coring{C}$, and
$\varepsilon_{\coring{C}} : \coring{C} \rightarrow A$ that satisfy
the equations:
\begin{equation}\label{coasociativa}
(\coring{C} \tensor{A} \Delta_{\coring{C}}) \circ
\Delta_{\coring{C}} = (\Delta_{\coring{C}} \tensor{A} \coring{C})
\circ \Delta_{\coring{C}}
\end{equation}
\begin{equation}\label{counitaria}
(\varepsilon_{\coring{C}} \tensor{A} \coring{C}) \circ
\Delta_{\coring{C}} = d_{\coring{C}}^-, \qquad (\coring{C}
\tensor{A} \varepsilon_{\coring{C}}) \circ \Delta_{\coring{C}} =
d_{\coring{C}}^+
\end{equation}
In \eqref{coasociativa} we have considered the canonical isomorphism
$\coring{C} \tensor{A} (\coring{C} \tensor{A} \coring{C}) \cong
(\coring{C} \tensor{A} \coring{C}) \tensor{A} \coring{C}$ as one
equality, denoting the ``common value'' by $\coring{C} \tensor{A}
\coring{C} \tensor{A} \coring{C}$. When $A$ is unital, we recover
the original definition from \cite{Sweedler:1975}. A direct
computation, using \eqref{coasociativa} and \eqref{counitaria},
shows that each $A$--coring
$(\coring{C},\Delta_{\coring{C}},\varepsilon_{\coring{C}})$
determines a comonad on $\rmod{A}$ defined by the functor
\begin{equation*}
\xymatrix{ - \tensor{A} \coring{C} : \rmod{A} \ar[r] & \rmod{A}}
\end{equation*}
and the natural transformations
\begin{equation*}
\xymatrix{- \tensor{A} \coring{C} \ar^-{-\tensor{A}
\Delta_{\coring{C}}} [rr] & & - \tensor{A} \coring{C} \tensor{A}
\coring{C} \\
- \tensor{A} \coring{C} \ar^-{- \tensor{A}
\varepsilon_{\coring{C}}}[rr] & & - \tensor{A} A \cong id }
\end{equation*}
The category of coalgebras for this comonad is the category
$\rcomod{\coring{C}}$ of all right $\coring{C}$--comodules.

Let $B$ be a firm ring, and a firm $B-A$--bimodule $\Sigma$.
Arguing as in \cite{Gomez/Vercruysse:2005arXiv}, we have a pair of
adjoint functors
\begin{equation}\label{adjcoring}
\xymatrix{\rmod{B} \ar@<0.5ex>^-{\dostensor{-}{B}{\Sigma}}[rrr] & &
& \rmod{A}, \ar@<0.5ex>^-{\dostensor{\hom{A}{\Sigma}{-}}{B}{B}}[lll]
& \dostensor{-}{B}{\Sigma} \dashv
\dostensor{\hom{A}{\Sigma}{-}}{B}{B}}
\end{equation}
which is obtained by composing the adjunctions
\begin{equation*}
\xymatrix{\rmodu{B} \ar@<0.5ex>^-{\dostensor{-}{B}{\Sigma}}[rrr] &
& & \rmod{A}, \ar@<0.5ex>^-{\hom{A}{\Sigma}{-}}[lll] &
\dostensor{-}{B}{\Sigma} \dashv \hom{A}{\Sigma}{-}}
\end{equation*}
and
\begin{equation*}
\xymatrix{\rmod{B} \ar@<0.5ex>^-{J}[rrr] & & & \rmodu{B},
\ar@<0.5ex>^-{\dostensor{-}{B}{B}}[lll] & J \dashv
\dostensor{-}{B}{B}}
\end{equation*}
where $J : \rmod{B} \rightarrow \rmodu{B}$ is the inclusion
functor.
  It will be useful
to give explicitly the unit and the counit of the adjunction
\eqref{adjcoring}. To do this, given $y \in Y$ for $Y$ a right
firm $B$--module, we will use the notation $d_Y^+(y) = y^b
\tensor{B} b \in Y \tensor{B} B$ (sum understood). Of course, this
element of the tensor product is determined by the condition $y^bb
= y$. The counit of the adjunction is
\begin{equation}\label{unit3}
\xymatrix{\eta_Y : Y \ar[r] &
\dostensor{\hom{A}{\Sigma}{\dostensor{Y}{B}{\Sigma}}}{B}{B}, &
\eta_Y(y) = \dostensor{(\dostensor{y^b}{B}{-})}{B}{b},}
\end{equation}
and the unit
\begin{equation}\label{counit3}
\xymatrix{\epsilon_X :
\trestensor{\hom{A}{\Sigma}{X}}{B}{B}{B}{\Sigma} \ar[r] & X, &
\epsilon_X(\trestensor{f}{B}{b}{B}{u}) = f(bu)}.
\end{equation}
We have then the associated comonad
\begin{equation*}
(\trestensor{\hom{A}{\Sigma}{-}}{B}{B}{B}{\Sigma},
\dostensor{\eta_{\dostensor{\hom{A}{\Sigma}{-}}{B}{B}}}{B}{\Sigma},\epsilon)
\end{equation*}
Assume we have a structure of $B - \coring{C}$--bicomodule over
$\Sigma$, that is, a homomorphism of $B-A$--bimodules
$\varrho_\Sigma : \Sigma \rightarrow \Sigma \tensor{A} \coring{C}$
such that
\begin{equation}\label{rcomod}
(\varrho_{\Sigma} \tensor{A} \coring{C}) \circ \varrho_{\Sigma} =
(\Sigma \tensor{A} \Delta_{\coring{C}}) \circ \varrho_{\Sigma},
\qquad (\Sigma \tensor{A} \varepsilon_{\coring{C}})\circ
\varrho_{\Sigma} = d_{\Sigma}^+
\end{equation}
It follows from \eqref{rcomod} that the natural transformation
\begin{equation}\label{(C)}
\xymatrix{\beta : - \tensor{B} \Sigma \ar^-{- \tensor{B}
\varrho_{\Sigma}}[rr] & & - \tensor{B} \Sigma \tensor{A}
\coring{C}}
\end{equation}
satisfies the conditions of the statement \emph{(C)} of the
Proposition \ref{correspondencia} and, therefore, it gives rise to
a canonical homomorphism of comonads $\mathsf{can}$ defined by the
composite
\begin{equation*}
\xymatrix{\trestensor{\hom{A}{\Sigma}{-}}{B}{B}{B}{\Sigma}
\ar^-{\trestensor{\hom{A}{\Sigma}{-}}{B}{B}{B}{\varrho_{\Sigma}}}
[rrr] \ar_-{\mathsf{can}}[rrrd] & & &
\fourtensor{\hom{A}{\Sigma}{-}}{B}{B}{B}{\Sigma}{A}{\coring{C}}
\ar^{\epsilon \tensor{A} \coring{C}}[d]
\\ & & & - \tensor{A} \coring{C},}
\end{equation*}
or, using a simplified version of Heynemann-Sweedler's notation, we
have, for each right $B$--module $X$:
\begin{equation*}
\xymatrix{\trestensor{\hom{A}{\Sigma}{X}}{B}{B}{B}{\Sigma}
\ar^-{\mathsf{can}_X}[rr] & & X \tensor{A} \coring{C}
\\
\trestensor{f}{B}{b}{B}{u} \ar@{|->}[rr] & &
f(bu_{[0]})\tensor{A}u_{[1]} }
\end{equation*}
where $\varrho_{\Sigma}(u) = u_{[0]} \tensor{A} u_{[1]}$ (sum
understood).

\medskip

 We can now apply Proposition \ref{adjuncion} and Theorem
\ref{Descent} to obtain:

\begin{theorem}\label{debil}
The functor $ - \tensor{B} \Sigma : \rmod{B} \rightarrow
\rcomod{\coring{C}}$ has a right adjoint
\begin{equation*}
\xymatrix{\dostensor{\hom{\coring{C}}{\Sigma}{-}}{B}{B} :
\rcomod{\coring{C}} \ar[rr] & &  \rmod{B}}
\end{equation*}
This functor is faithful and full if and only if $\mathsf{can}$ is
a natural isomorphism and $- \tensor{B} \Sigma : \rmod{B}
\rightarrow \rmod{A}$ preserves the equalizer
\begin{equation}\label{iguamodu}
\xymatrix{ \hom{\coring{C}}{\Sigma}{X} \tensor{B} B \ar[r]  &
\hom{A}{\Sigma}{X}\tensor{B} B \ar@<.5ex>^-{\alpha_X}[rrr]
\ar@<-.5ex>_-{\hom{A}{\Sigma}{\varrho_X}\tensor{B}{B}}[rrr] & & &
\hom{A}{\Sigma}{X \tensor{A} \coring{C}} \tensor{B} B}
\end{equation}
for every right $\coring{C}$--comodule $(X,\varrho_X)$, where
$\alpha_X(f \tensor{B} b) = [(f \tensor{A} \coring{C}) \circ
\varrho_{\Sigma}] \tensor{B} b$ for $f \tensor{B} b \in
\hom{A}{\Sigma}{X} \tensor{B} B$.
\end{theorem}
\begin{proof}
It is enough to check that, in the present situation, $\alpha_X$,
as defined in \eqref{rhophi}, is given as in the statement of the
theorem, and that $- \tensor{B} B : \rmodu{B} \rightarrow
\rmod{B}$ is left exact since it is right adjoint to the inclusion
functor $J : \rmod{B} \rightarrow \rmodu{B}$.
\end{proof}

Theorem \ref{comonadic} has the following consequence:

\begin{theorem}\label{fuerte}
Let $\Sigma$ be a $B-\coring{C}$--bicomodule. The functor $ -
\tensor{B} \Sigma : \rmod{B} \rightarrow \rcomod{\coring{C}}$ is
an equivalence of categories if and only if $- \tensor{B} \Sigma :
\rmod{B} \rightarrow \rmod{A}$ preserves the equalizers of the
form \eqref{iguamodu}, reflects isomorphisms, and the canonical
transformation $\mathsf{can}$ is an isomorphism.
\end{theorem}

Every right $\coring{C}$--comodule $\Sigma$, is a
$T-\coring{C}$--bicomodule, where $T = \rend{\coring{C}}{\Sigma}$.
Of course, in this case we have a natural isomorphism
$\hom{A}{\Sigma}{-}\tensor{T} T \cong \hom{A}{\Sigma}{-}$. Following
\cite{Wisbauer:2004arXiv}, we will say that $\Sigma$ is a
\emph{Galois $\coring{C}$--comodule} when $ \mathsf{can} :
\hom{A}{\Sigma}{-} \tensor{T} \Sigma \rightarrow - \tensor{A}
\coring{C} $ is an isomorphism. More generally, if $\Sigma$ is a
$B-\coring{C}$--bicomodule firm as a $B-A$--bimodule, then we have a
homomorphism of rings $\lambda : B \rightarrow T$. Assume that $B$
is a left ideal of $T$ (that is, that the image under $\lambda$ of
$B$ is a left ideal or $T$). By \cite[Lemma 4.10, Lemma
4.11]{Gomez/Vercruysse:2005arXiv}, we have a commutative diagram of
natural transformations:
\begin{equation}\label{cancan}
\xymatrix{\hom{A}{\Sigma}{-} \tensor{B} B \tensor{B} \Sigma
\ar^{\simeq}[rr] \ar[rd]_{\mathsf{can}}& &  \hom{A}{\Sigma}{-}\tensor{T} \Sigma \ar[ld]^{\mathsf{can}} \\
& - \tensor{A} \coring{C} &}
\end{equation}
Therefore, $\Sigma$ is a Galois right $\coring{C}$--comodule in the
sense of \cite{Wisbauer:2004arXiv} if and only if
\begin{equation*}
\mathsf{can} : \hom{A}{\Sigma}{-} \tensor{B} B \tensor{B} \Sigma
\rightarrow - \tensor{A} \coring{C}
\end{equation*}
 is an isomorphism. Obviously, Theorems
\ref{debil} and \ref{fuerte} can be formulated in the case where $B
= T = \rend{\coring{C}}{\Sigma}$, characterizing then when
$\hom{\coring{C}}{\Sigma}{-} : \rcomod{\coring{C}} \rightarrow
\rmod{T}$ is full and faithful or it gives an equivalence of
categories. This evokes Gabriel-Popescu's Theorem
\cite{Gabriel/Popescu:1964}.

A left module $M$ over a firm ring $R$ will be said to be
\emph{flat} if the functor $ - \tensor{R} M : \rmodu{R} \rightarrow
\mathsf{Ab}$ is exact, where $\mathsf{Ab}$ denotes the category of
abelian groups.

The following theorem gives a general version, in the sense that no
finiteness condition is assumed a priori on the comodule $\Sigma$,
of \cite[Theorem 2.1.(1)]{Brzezinski:2005}, \cite[Theorem
3.8]{Caenepeel/DeGroot/Vercruysse:unp2005},
\cite[18.27]{Brzezinski/Wisbauer:2003}, \cite[Theorem
4.9]{Gomez/Vercruysse:2005arXiv}.

\begin{theorem}\label{fielmenteplano}
Let $\Sigma$ be a $B-\coring{C}$--bicomodule, where $A$ is unital
and $\coring{C}$ is an $A$--coring, flat as a left $A$--module.
Assume $B$ to be firm, and that ${}_B\Sigma$ is a firm module. If
$B$ is a left ideal of $T = \rend{\coring{C}}{\Sigma}$, then the
following statements are equivalent
\begin{enumerate}[(i)]
\item\label{GP1} The functor $\hom{\coring{C}}{\Sigma}{-} :
\rcomod{\coring{C}} \rightarrow \rmod{T}$ is full and faithful;
\item\label{GP2} $\Sigma$ is a generator of the category
$\rcomod{\coring{C}}$; \item\label{GP3} $\mathsf{can} :
\hom{\coring{C}}{\Sigma}{-} \tensor{T} \Sigma \rightarrow -
\tensor{A} \coring{C}$ is an isomorphism, and ${}_T\Sigma$ is flat;
\item\label{GP4} $\mathsf{can} : \hom{A}{\Sigma}{-} \tensor{B} B
\tensor{B} \Sigma \rightarrow - \tensor{A} \coring{C}$ is an
isomorphism of comonads on $\rmod{A}$ and ${}_B\Sigma$ is flat;
\item\label{GP5} the functor $\hom{\coring{C}}{\Sigma}{-} \tensor{B}
B : \rcomod{\coring{C}} \rightarrow \rmod{B}$ is full and faithful.
\end{enumerate}
\end{theorem}
\begin{proof}
By \cite[Proposition 1.2]{ElKaoutit/Gomez/Lobillo:2004},
${}_A\coring{C}$ is flat if and only if $\rcomod{\coring{C}}$ is a
Grothendieck category and the forgetful functor $U :
\rcomod{\coring{C}} \rightarrow
\rmod{A}$ is exact. \\
\eqref{GP1} $\Leftrightarrow$ \eqref{GP2} is a consequence of
Gabriel-Popescu's Theorem \cite[Cap. III, Teorem\u a
9.1.(2)]{Nastasescu:1976}.\\
\eqref{GP2} $\Rightarrow$ \eqref{GP3} By Theorem \ref{debil}, with
$B = T$, $\Sigma$ is a Galois $\coring{C}$--comodule. By
Gabriel-Popescu's Theorem \cite[Cap. III, Teorem\u a
9.1.(3)]{Nastasescu:1976} the functor $ - \tensor{T} \Sigma :
\rmod{T} \rightarrow \rcomod{\coring{C}}$ is exact. Since $U :
\rcomod{\coring{C}} \rightarrow \rmod{A}$ is exact, we get that
${}_T\Sigma$ is flat. \\
\eqref{GP3} $\Rightarrow$ \eqref{GP1}
By Theorem \ref{debil} with $B = T$. \\
\eqref{GP3} $\Rightarrow$ \eqref{GP4} In view of diagram
\eqref{cancan}, $\mathsf{can} : \hom{A}{\Sigma}{-} \tensor{B} B
\tensor{B} \Sigma \rightarrow - \tensor{A} \coring{C}$ is an
isomorphism. We have proved that \eqref{GP3} is equivalent to
\eqref{GP2}. Therefore, we deduce from \cite[Proposition
4.13]{Gomez/Vercruysse:2005arXiv}
that ${}_B\Sigma$ is flat. \\
\eqref{GP4} $\Rightarrow$ \eqref{GP5} By Theorem \ref{debil}.\\
\eqref{GP5} $\Rightarrow$ \eqref{GP2} Follows by the standard
argument: For any right $\coring{C}$--comodule $X$, take a free
presentation $T^{(I)} \rightarrow \hom{\coring{C}}{\Sigma}{X}
\rightarrow 0$. After tensoring on the right by $B \tensor{B}
\Sigma$, and taking that $T \tensor{B} B \cong B$ \cite[Lemma
4.10]{Gomez/Vercruysse:2005arXiv} into account, we get a
presentation $\Sigma^{(I)} \rightarrow X \rightarrow 0$.
\end{proof}

\begin{lemma}\label{sobreideal}
Let $\Sigma$ be a $B$-$\coring{C}$--bicomodule, where $B$ is a
firm ring. Assume that $\Sigma$ is firm as a left $B$--module. If
the counit of the adjunction $- \tensor{B} \Sigma \dashv
\hom{\coring{C}}{\Sigma}{-} \tensor{B} B$ evaluated at $B$ is
surjective, then $B$ is a left ideal of $T =
\rend{\coring{C}}{\Sigma}$.
\end{lemma}
\begin{proof}
The proof of the implication $(iv) \Rightarrow (v)$ of
\cite[Theorem 4.15]{Gomez/Vercruysse:2005arXiv} runs here: In view
of \eqref{unit3} and Remark \ref{unitcounit}, the hypothesis says
that the map
\begin{equation*}
\xymatrix{\eta_B : B \ar[r] & \hom{\coring{C}}{\Sigma}{B
\tensor{B} \Sigma} \tensor{B} B, & \eta_B(b) =
\dostensor{(\dostensor{b^c}{B}{-})}{B}{c}}
\end{equation*}
is surjective. By composing $\eta_B$ with the isomorphism
$\hom{\coring{C}}{\Sigma}{B \tensor{B} \Sigma} \cong
\hom{\coring{C}}{\Sigma}{\Sigma} \tensor{B} B = T \tensor{B} B$,
we obtain that every element of $T \tensor{B} B$ is a finite sum
of elements of the form $b^c \cdot \tensor{B} c$. Thus, $B$ is a
left ideal of $T$.
\end{proof}

\begin{theorem}\label{GE}
Let $\Sigma$ be a $B-\coring{C}$--bicomodule, where $A$ is unital
and $\coring{C}$ is an $A$--coring, flat as a left $A$--module.
Assume $B$ to be firm, and that ${}_B\Sigma$ is a firm module. The
following statements are equivalent:
\begin{enumerate}[(i)]
\item\label{GE1} The functor $- \tensor{B} \Sigma : \rmod{B}
\rightarrow
\rcomod{\coring{C}}$ is an equivalence of categories; %
\item\label{GE2} $\mathsf{can} : \hom{A}{\Sigma}{-} \tensor{B} B
\tensor{B} \Sigma \rightarrow - \tensor{A} \coring{C}$ is an
isomorphism, and ${}_B
\Sigma$ is faithfully flat; %
\item\label{GE3} $\Sigma$ is a generator of $\rcomod{\coring{C}}$
such that the functor $- \tensor{B} \Sigma : \rmod{B} \rightarrow
\rcomod{\coring{C}}$ is full and faithful; %
\item\label{GE4} $\Sigma$ is a generator of $\rcomod{\coring{C}}$
such that the functor $- \tensor{B} \Sigma : \rmod{B} \rightarrow
\rcomod{\coring{C}}$ is faithful and $B$ is a left ideal of $T$.
\end{enumerate}
\end{theorem}
\begin{proof}
\eqref{GE1} $\Rightarrow$ \eqref{GE2} $\mathsf{can}$ is an
isomorphism by Theorem \ref{fuerte}. By \cite[Proposition
1.2]{ElKaoutit/Gomez/Lobillo:2004}, the forgetful functor
$\rcomod{\coring{C}} \rightarrow \rmod{A}$ is faithful and exact.
Therefore, the functor $ - \tensor{B} \Sigma : \rmod{B}
\rightarrow \rmod{A}$
is faithful and exact. \\
\eqref{GE2} $\Rightarrow$ \eqref{GE1} By Theorem \ref{fuerte}.\\
\eqref{GE1} $\Rightarrow$ \eqref{GE3} Since $B$ is a generator of
$\rmod{B}$, $\Sigma \cong B \tensor{B} \Sigma$ is a generator of
$\rcomod{\coring{C}}$.\\
\eqref{GE3} $\Rightarrow$ \eqref{GE4} By Lemma \ref{sobreideal}.\\
\eqref{GE4} $\Rightarrow$ \eqref{GE2} The forgetful functor
$\rcomod{\coring{C}} \rightarrow \rmod{A}$ is faithful. Therefore,
$- \tensor{B} \Sigma : \rmod{B} \rightarrow \rmod{A}$ is faithful.
By Theorem \ref{fielmenteplano}, ${}_B\Sigma$ is flat and
$\mathsf{can}$ is an isomorphism.
\end{proof}

\begin{remark}
Theorem \ref{GE} shows, in conjunction with Theorem
\ref{fielmenteplano}, that if $- \tensor{B} \Sigma : \rmod{B}
\rightarrow \rcomod{\coring{C}}$ is an equivalence of categories,
then $\Sigma$ is a Galois comodule in the sense of
\cite{Wisbauer:2004arXiv}.
\end{remark}

In the particular case where $B$ is unital and $\Sigma_A$ is
finitely generated and projective, we have a natural isomorphism
\begin{equation*}
\nu : \hom{A}{\Sigma}{-} \tensor{B} B \cong - \tensor{A} \Sigma^*
\qquad ( f \tensor{B} b \mapsto f(be_i) \tensor{A} e_i^*)
\end{equation*}
where $\Sigma^* = \hom{A}{\Sigma}{A}$, and $\{ (e_i, e_i^*) \}$ is a
finite dual basis for $\Sigma_A$. We have then the adjoint pair $-
\tensor{B} \Sigma \dashv - \tensor{A} \Sigma^*$ and the associated
comonad $- \tensor{A} \Sigma^* \tensor{B} \Sigma$ on $\rmod{A}$. In
particular, $\Sigma^* \tensor{B} \Sigma$ is an $A$--coring, the
\emph{comatrix coring} associated to the bimodule ${}_B\Sigma_A$
(see \cite{ElKaoutit/Gomez:2003}). If $\Sigma$ is a
$B$-$\coring{C}$--bicomodule, then we have a commutative diagram of
homomorphisms of comonads
\begin{equation}\label{cancan2}
\xymatrix{\hom{A}{\Sigma}{-} \tensor{B} B \tensor{B} \Sigma \ar^{\nu
\tensor{B} \Sigma}[rr]
\ar_{\mathsf{can}}[dr] && - \tensor{A} \Sigma^* \tensor{B} \Sigma \ar^{- \tensor{A} \mathbf{can}}[dl]\\
&- \tensor{A} \coring{C},&}
\end{equation}
where $\mathbf{can} (\phi \tensor{B} u) = \phi(u_{[0]})u_{[1]}$ for
$\phi \tensor{B} u \in \Sigma^* \tensor{B} \Sigma$ is a homomorphism
of $A$--corings from $\Sigma^* \tensor{B} \Sigma$ to $\coring{C}$
(see \cite{ElKaoutit/Gomez:2003}).

\begin{corollary}\cite[Theorem 3.2, Theorem
3.10]{ElKaoutit/Gomez:2003}\label{clasico} Let $\Sigma$ be a
$B-\coring{C}$--bicomodule, where $B$ is a unital ring and
$\coring{C}$ is a coring over a unital ring $A$. The following
statements are equivalent:
\begin{enumerate}[(i)]
\item\label{PE1} ${}_A\coring{C}$ is flat and the functor $- \tensor{B} \Sigma
: \rmod{B} \rightarrow \rcomod{\coring{C}}$ is an equivalence of
categories;
\item\label{PE2} $\Sigma_A$ is finitely generated and projective, the canonical map $\mathbf{can} : \Sigma^*
\tensor{B} \Sigma \rightarrow \coring{C}$ is an isomorphism, and
${}_B\Sigma$ is faithfully flat;
\item\label{PE3} ${}_A\coring{C}$ is flat, $\Sigma$ is a finitely generated
projective generator of $\rcomod{\coring{C}}$, and $\lambda : B
\rightarrow T$ is an isomorphism.
\end{enumerate}
\end{corollary}
\begin{proof}
\eqref{PE1} $\Rightarrow$ \eqref{PE3} The equivalence of categories
$- \tensor{B} \Sigma$ preserves finitely generated and projective
generators, whence $\Sigma \cong B \tensor{B} \Sigma$ is a finitely
generated projective generator of $\rcomod{\coring{C}}$. On the
other hand, by Theorem \ref{GE}, $B$ is a left ideal of $T$. Since
$B$ is a unital subring of $T$, we must have $B = T$. \\
\eqref{PE3} $\Rightarrow$ \eqref{PE2} Since the forgetful functor $U
: \rcomod{\coring{C}} \rightarrow \rmod{A}$ has an exact right
adjoint $- \tensor{A} \coring{C}$ which preserves colimits, $U$
preserves finitely generated and projective objects. The rest
follows directly from Theorem
\ref{GE} and diagram \eqref{cancan2}. \\
\eqref{PE2} $\Rightarrow$ \eqref{PE1} Since $\mathbf{can}$ is an
isomorphism, and the left $A$--module $\Sigma^* \tensor{B} \Sigma$
is flat, we get that ${}_A\coring{C}$ is flat. That $ - \tensor{B}
\Sigma$ is an equivalence of categories follows from Theorem
\ref{GE} and diagram \eqref{cancan2}.
\end{proof}

\begin{remark}
It follows from Corollary \ref{clasico} that if $- \tensor{B} \Sigma
: \rmod{B} \rightarrow \rcomod{\coring{C}}$ is an equivalence, and
$B$ is unital, then the coring $\coring{C}$ is Galois in the sense
of \cite{ElKaoutit/Gomez:2003} (or $\Sigma$ is a Galois comodule as
in \cite{Brzezinski/Wisbauer:2003}), that is, $\mathbf{can} :
\Sigma^* \tensor{T} \Sigma \rightarrow \coring{C}$ is an
isomorphism.
\end{remark}

Next, we will derive some fundamental results of
\cite{Gomez/Vercruysse:2005arXiv}. Write $\Sigma^* =
\hom{A}{\Sigma}{A}$, and assume that $A$ is a ring with unit. The
$B$--bimodule $S = \Sigma \tensor{A} \Sigma^*$ has the structure
of a $B$--ring (without unit, in general), with associative
multiplication given by
\begin{equation*}
\mu(x \tensor{A} \phi \tensor{B} y \tensor{A} \psi) = x \phi(y)
\tensor{A} \psi = x \tensor{A} \phi(y) \psi, \qquad (x \tensor{A}
\phi \tensor{B} y \tensor{A} \psi \in \Sigma \tensor{A} \Sigma^*
\tensor{B} \Sigma \tensor{A} \Sigma^*)
\end{equation*}
The map $S \rightarrow \rend{A}{\Sigma}$ sending
$\dostensor{x}{A}{\phi}$ onto the endomorphism $y \mapsto
x\phi(y)$ is a homomorphism of $B$--rings. In particular, $\Sigma$
is an $S-A$--bimodule. Analogously, the homomorphism $S
\rightarrow \lend{A}{\Sigma^*}$ that maps $\dostensor{y}{A}{\psi}$
to the endomorphism $\phi \mapsto \phi(y)\psi$ is a homomorphism
of $B$--rings, and $\Sigma^*$ is then an $A-S$--bimodule.
 The situation studied in
\cite{Gomez/Vercruysse:2005arXiv} begins with a homomorphism of
rings $\iota : R \rightarrow S$, where $R$ is a firm ring.  We will
use the notation $\iota (r) = e_r \tensor{R} e^*_r$ for $r \in R$
(sum understood). By restriction of scalars, $\Sigma$ is then an
$R-A$--bimodule, with the left action of $R$ defined explicitly by
$ru = e_re_r^*(u)$, for $r \in R$, $u \in \Sigma$. The
$A-R$--bimodule structure on $\Sigma^*$ is given by the right action
of $R$ defined by $\phi r = \phi (e_r)e_r^*$. Now, $S = \Sigma
\tensor{A} \Sigma^*$ becomes an $R$--ring, and thus we could take $B
= R$ without loss of generality (this was not formally the case in
\cite{Gomez/Vercruysse:2005arXiv}, where the initial ring $B$ was
assumed to be unital).

If $\Sigma$ is firm as a left $R$--module, then we have, following
\cite{Gomez/Vercruysse:2005arXiv}, a natural isomorphism
\begin{equation}\label{isomor}
\nu : \hom{A}{\Sigma}{-} \tensor{R} R \simeq - \tensor{A} \Sigma^*
\tensor{R} R \qquad (h \tensor{R} r \mapsto h(e_s) \tensor{A} e^*_s
\tensor{B} r^s)
\end{equation}

In view of \eqref{isomor} the adjunction \eqref{adjcoring} leads,
using the notation $\Sigma^{\dag} = \Sigma^* \tensor{R} R$, to an
adjunction
\begin{equation}\label{adj6}
\xymatrix{\rmod{R} \ar@<0.5ex>^-{\dostensor{-}{R}{\Sigma}}[rrr] &
& & \rmod{A}, \ar@<0.5ex>^-{\dostensor{-}{A}{\Sigma^{\dag}}}[lll]
& \dostensor{-}{R}{\Sigma} \dashv \dostensor{-}{A}{\Sigma^{\dag}}}
\end{equation}
whose counit is
\begin{equation}\label{counit6}
\xymatrix{\epsilon_M : \trestensor{M}{A}{\Sigma^{\dag}}{R}{\Sigma}
\ar[r] & M & (\fourtensor{m}{A}{\phi}{R}{r}{R}{x} \mapsto
m\phi(rx)),}
\end{equation}
and whose unit is
\[
\xymatrix{\eta_N : N \ar[r] &
\trestensor{N}{R}{\Sigma}{A}{\Sigma^{\dag}} & (n \mapsto
\fourtensor{n^r}{R}{e_s}{A}{e^*_s}{R}{r^s})}
\]
From the adjoint pair \eqref{adj6} we have the comonad on
$\rmod{A}$:
\begin{equation}\label{comatrixcomon}
(- \tensor{A} \Sigma^{\dag} \tensor{R} \Sigma, \eta_{- \tensor{A}
\Sigma^{\dag}} \tensor{R} \Sigma, \epsilon )
\end{equation}
Evaluating at $A$, we obtain the $A$--coring $A \tensor{A}
\Sigma^{\dag} \tensor{R} \Sigma \cong \Sigma^{\dag} \tensor{R}
\Sigma$ with comultiplication $\Delta^{\dag} = \eta_{A \tensor{A}
\Sigma^{\dag}} \tensor{R} \Sigma$ and counity $\epsilon_A$.
Explicitly,
\begin{equation*}
\Delta^{\dag}(\trestensor{\phi}{R}{r}{R}{u}) =
\fourtensor{\phi}{R}{s}{R}{\abrir{t}}{R}\dostensor{(r^s)^t}{R}{u}
\end{equation*}
\begin{equation*}
\epsilon_A(\phi \tensor{R} r \tensor{B} u) = \phi(ru)
\end{equation*}
Moreover, the comonad \eqref{comatrixcomon} is determined by the
\emph{comatrix coring} $(\Sigma^{\dag} \tensor{R} \Sigma,
\Delta^{\dag},\epsilon_A)$, since the functor underlying to the
comonad is a tensor product over $A$.

Returning to the case where $(\Sigma, \varrho_{\Sigma})$ is an
$R-\coring{C}$--bicomodule, we have that the natural transformation
\eqref{(C)} gives rise, by Proposition \ref{correspondencia}, to a
homomorphism of comonads $\mathsf{can}^{\dag}$ defined (see
\eqref{canr}) at $(X,\varrho_X)$ by
\begin{equation*}
{\mathsf{can}}^{\dag}_X = (X \tensor{A} \epsilon_A \tensor{A}
\coring{C}) \circ (X \tensor{A} \Sigma^{\dag} \tensor{R}
\varrho_{\Sigma})
\end{equation*}
Obviously,
\begin{equation}\label{canplus}
\mathsf{can}^{\dag}_X = X \tensor{A} \mathbf{can}^{\dag}
\end{equation}
where $\mathbf{can}^{\dag}$ is the map
\begin{equation}
\xymatrix{\mathbf{can}^{\dag} : \Sigma^{\dag} \tensor{R} \Sigma
\ar[r] & \coring{C}, & & \phi \tensor{B} r \tensor{R} u \mapsto
\phi(ru_{[0]})u_{[1]}}
\end{equation}
which turns out to be a homomorphism of $A$--corings, because
$\mathsf{can}^{\dag}$ is a homomorphism of comonads. According to
Proposition \ref{adjuncion}, the functor $- \tensor{R} \Sigma :
\rmod{R} \rightarrow \rcomod{\coring{C}}$ has a right adjoint that,
over a right $\coring{C}$--comodule $(X, \varrho_X)$, is defined as
the equalizer in $\rmod{R}$ of the pair $(\alpha_X,\varrho_X
\tensor{A} \Sigma^{\dag})$. An easy computation shows that, in the
present case, $\alpha_X = X \tensor{A} \alpha_{\Sigma^{\dag}}$,
where
\begin{equation}\label{sigmaplus}
\alpha_{\Sigma^{\dag}} = (\mathbf{can}^{\dag} \tensor{A}
\Sigma^{\dag}) \circ \eta_{\Sigma^{\dag}}
\end{equation}
From the diagrams \eqref{comod} we get that
$(\Sigma^{\dag},\alpha_{\Sigma^{\dag}})$ is a
$\coring{C}-R$--bicomodule. In this way, the functor defined in
Proposition \ref{adjuncion} becomes a cotensor product, as it is
defined by the equalizer
\begin{equation}\label{cotensor}
\xymatrix{X \cotensor{\coring{C}} \Sigma^{\dag} \ar[r] & X
\tensor{A} \Sigma^{\dag} \ar@<.5ex>^-{X \tensor{A}
\alpha_{\Sigma^{\dag}}}[rr] \ar@<-.5ex>_-{\varrho_X \tensor{A}
\Sigma^{\dag}}[rr] & & X \tensor{A} \coring{C} \tensor{A}
\Sigma^{\dag} }
\end{equation}
From Proposition \ref{adjuncion} and Theorem \ref{Descent} we
deduce:

\begin{theorem}\cite[Theorem 4.9]{Gomez/Vercruysse:2005arXiv} Let $\coring{C}$
be an $A$--coring and $\Sigma$ a right $\coring{C}$--comodule, which
is as well a left $B$--module. Let $\iota : R \rightarrow \Sigma
\tensor{A} \Sigma^*$ be a homomorphism of rings, where $R$ is a firm
ring. If $\Sigma$ is an $R-\coring{C}$--bicomodule such that
$\Sigma$ is firm as a left $R$--module, then the functor $-
\tensor{R} \Sigma : \rmod{R} \rightarrow \rcomod{\coring{C}}$ admits
as a right adjoint the cotensor product functor $-
\cotensor{\coring{C}} \Sigma^{\dag} : \rcomod{\coring{C}}
\rightarrow \rmod{R}$. This functor is faithful and full if and only
if $\mathbf{can}^{\dag} : \Sigma^{\dag} \tensor{R} \Sigma
\rightarrow \coring{C}$ is an isomorphism of $A$--corings and $ -
\tensor{R} \Sigma : \rmod{R} \rightarrow \rmod{A}$ preserves the
equalizers of the form \eqref{cotensor}.
\end{theorem}

Theorem \ref{comonadic} boils down to:

\begin{theorem}  Let $\coring{C}$
be an $A$--coring and $\Sigma$ a right $\coring{C}$--comodule, which
is as well a left $B$--module. Let $\iota : R \rightarrow \Sigma
\tensor{A} \Sigma^*$ be a homomorphism of rings, where $R$ is a firm
ring. If $\Sigma$ is an $R-\coring{C}$--bicomodule such that
$\Sigma$ is firm as a left $R$--module, then the functor $-
\tensor{R} \Sigma : \rmod{R} \rightarrow \rcomod{\coring{C}}$ is an
equivalence of categories if and only if $\mathbf{can}^{\dag} :
\Sigma^{\dag} \tensor{R} \Sigma \rightarrow \coring{C}$ is an
isomorphism of $A$--corings, and $ - \tensor{R} \Sigma : \rmod{R}
\rightarrow \rmod{A}$ reflects isomorphisms and preserves the
equalizers of the form \eqref{cotensor}.
\end{theorem}

\begin{remark}
Start from a $B-\coring{C}$--bicomodule $\Sigma$, with $B$ and
${}_B\Sigma$ firm, and consider a ring homomorphism $\iota : B
\rightarrow S = \Sigma \tensor{A} \Sigma^*$ (that is, put $B = R$ in
the foregoing discussion, which essentially does not constitute a
restriction, as we discussed before). This is just to say that the
left multiplication ring homomorphism $\lambda : B \rightarrow T =
\rend{\coring{C}}{\Sigma}$ factorizes throughout $S$, in the sense
of the commutative diagram
\begin{equation*}
\xymatrix{\Sigma \tensor{A} \Sigma^* \ar[rr] & & \rend{A}{\Sigma}
\\
B \ar^{\iota}[u] \ar^{\lambda}[rr] & & \rend{\coring{C}}{\Sigma}
\ar[u]}
\end{equation*}
Now, since in this case the functors $\hom{\coring{C}}{\Sigma}{-}
\tensor{B} B, - \cotensor{\coring{C}} \Sigma^{\dag} :
\rcomod{\coring{C}} \rightarrow \rmod{B}$ are both right adjoint to
$ - \tensor{B} \Sigma : \rmod{B} \rightarrow \rcomod{\coring{C}}$,
then they are isomorphic. On the other hand, from the natural
isomorphism $\nu$ given in \eqref{isomor}, we get the commutative
diagram of homomorphisms of comonads
\begin{equation}\label{cancan3}
\xymatrix{\hom{A}{\Sigma}{-} \tensor{B} B \tensor{B} \Sigma \ar^{\nu
\tensor{B} \Sigma}[rr]
\ar_{\mathsf{can}}[dr] && - \tensor{A} \Sigma^{\dag} \tensor{B} \Sigma \ar^{- \tensor{A} \mathbf{can}^{\dag}}[dl]\\
&- \tensor{A} \coring{C},&}
\end{equation}
Therefore, $\mathsf{can}$ is a natural isomorphism if and only if
$\mathbf{can}$ is an isomorphism \cite[Theorem
4.2]{Vercruysse:2006arXiv}. In this way, we can add the following
two statements to the list of equivalent conditions in Theorem
\ref{fielmenteplano}, connecting Theorem \ref{fielmenteplano} and
\cite[Theorem 4.9]{Gomez/Vercruysse:2005arXiv}:

\medskip

\noindent $(iv')$ \textit{$\mathbf{can} : \Sigma^{\dag} \tensor{B}
\Sigma \rightarrow \coring{C}$ is an isomorphism of corings and
${}_B\Sigma$ is flat;}\\
$(v')$ \textit{the functor $ - \cotensor{\coring{C}} \Sigma^{\dag} :
\rcomod{\coring{C}} \rightarrow \rmod{B}$ is full and faithful.}

\medskip

Analogously, in Theorem \ref{GE} we can add the following equivalent
condition, connecting Theorem \ref{GE} and \cite[Theorem
4.15]{Gomez/Vercruysse:2005arXiv}:

\medskip
\noindent
 $(ii')$ \textit{$\mathbf{can} : \Sigma^{\dag} \tensor{B}
\Sigma \rightarrow \coring{C}$ is an isomorphism of corings and
${}_B\Sigma$ is faithfully flat.}
\end{remark}

\section*{Acknowledgements}
I wish to thank Claudia Menini and Joost Vercruysse for a careful
reading of the (Spanish) earlier versions of this paper. Their
comments, in conjunction with that of the referee, contributed to
the improvement of this final (English) version. I am also grateful
to Bachuki Mesablishvili for drawing my attention to Dubuc's
monograph \cite{Dubuc:1970}, and to the referee for a question
concerning the notion of a Galois comodule from
\cite{Wisbauer:2004arXiv}, whose answer led me to formulate Theorems
\ref{fielmenteplano} and \ref{GE}.

\providecommand{\bysame}{\leavevmode\hbox
to3em{\hrulefill}\thinspace}
\providecommand{\MR}{\relax\ifhmode\unskip\space\fi MR }
\providecommand{\MRhref}[2]{%
  \href{http://www.ams.org/mathscinet-getitem?mr=#1}{#2}
} \providecommand{\href}[2]{#2}

\end{document}